\documentclass{amsart}
\usepackage{lastpage}
\usepackage{yhmath,verbatim,mathabx}
\usepackage[all]{xy} 

\renewcommand{\overset}[2]{%
	\mathop{#2}\limits^{\vbox to +.2ex{%
			\kern -1.ex\hbox{$\scriptstyle #1$}\vss}}}

\DeclareFontFamily{U}{mathx}{\hyphenchar\font45}
\DeclareFontShape{U}{mathx}{m}{n}{
	<5> <6> <7> <8> <9> <10>
	<10.95> <12> <14.4> <17.28> <20.74> <24.88>
	mathx10
}{}
\DeclareSymbolFont{mathx}{U}{mathx}{m}{n}
\DeclareFontSubstitution{U}{mathx}{m}{n}
\DeclareMathAccent{\widecheck}{0}{mathx}{"71}
\DeclareMathAccent{\wideparen}{0}{mathx}{"75}

\newcommand{\bdis}{\begin{displaymath}}
\newcommand{\edis}{\end{displaymath}}
\newcommand{\be}{\begin{equation}}
\newcommand{\ee}{\end{equation}}
\newcommand{\mbb}{\mathbb}
\newcommand{\mcal}{\mathcal}

\newcommand{\vp}{\varphi}

\newcommand{\zf}{\zeta\left(\frac{1}{2}+it\right)}

\DeclareMathOperator{\cn}{cn}

\theoremstyle{definition}

\theoremstyle{remark}
\newtheorem{remark}[]{Remark}

\newtheorem*{mydef11}{{\bf Theorem 1}}

\newtheorem*{mydef12}{{\bf Theorem 2}}

\newtheorem*{mydef13}{{\bf Theorem 3}}

\newtheorem*{mydef41}{{\bf Corollary 1}}

\newtheorem*{mydef42}{{\bf Corollary 2}}

\newtheorem*{mydef51}{{\bf Lemma 1}}

\newtheorem*{mydef52}{{\bf Lemma 2}}

\newtheorem*{mydef53}{{\bf Lemma 3}}

\numberwithin{equation}{section}

\begin{document}

\title{Jacob's ladders, crossbreeding and infinite sets of meta-functional equations as new species generated by the mother formula} 

\author{Jan Moser}

\address{Department of Mathematical Analysis and Numerical Mathematics, Comenius University, Mlynska Dolina M105, 842 48 Bratislava, SLOVAKIA}

\email{jan.mozer@fmph.uniba.sk}

\keywords{Riemann zeta-function}

\begin{abstract}
In this paper we obtain a set of meta-functional equations as new species of formulas in classical mathematical analysis. Mentioned species are generated by crossbreeding complete hybrid formula as a mother formula. Namely, they are generated by an infinite set of crossbreedings on some subsidiary infinite set of meta-functional equations with one neutral factor. \\ 

\begin{center} 
	DEDICATED TO THE 160th ANNIVERSARY OF DARWIN'S ORIGIN OF SPECIES	
\end{center} 
\end{abstract}
\maketitle

\section{Introduction} 

\subsection{} 

Let us remind that the following sets of values 
\be \label{1.1} 
\begin{split}
& \left\{ \left|\zf\right|^2 \right\},\ \{f_1(t)\}=\{\sin^2t\},\ \{f_2(t)\}=\{\cos^2t\}, \\ 
& t\in [\pi L,\pi L+U],\ U\in (0,\pi/2),\ L\in\mbb{N} 
\end{split}
\ee  
generate the exact complete hybrid formula (see \cite{7}, (3.2), $k_1=k_2=1$) 
\be \label{1.2} 
\tilde{Z}^2(\alpha_1^{1,1})\sin^2\alpha_0^{1,1}+\tilde{Z}^2(\alpha_0^{2,1})\cos^2\alpha_0^{2,1}=\tilde{Z}^2(\beta_1^1), 
\ee  
where  (comp. \cite{9}, (1.6)--(1.10)) 
\be \label{1.3} 
\begin{split}
& \alpha_r^{l,1}=\alpha_r(U,\pi L;f_l,|\zeta_{0,5}|^2),\ r=0,1,\ l=1,2, \\ 
& |\zeta_{0,5}|^2=\left|\zf\right|^2, \\ 
& \beta_1^1=\beta_1(U,\pi L;|\zeta_{0,5}|^2), \\
& \alpha_0^{l,1}\in (\pi L,\pi L+U),\ \alpha_1^{l,1},\beta_1^1\in (\overset{1}{\wideparen{\pi L}},\overset{1}{\wideparen{\pi L+U}}), \\ 
& U\in (0,\pi/2),\ L\geq L_0>0, 
\end{split}
\ee  
($L_0$ is sufficiently big one), next, we denote by the symbol 
\bdis 
[\overset{1}{\wideparen{\pi L}},\overset{1}{\wideparen{\pi L+U}}]
\edis  
the first reverse iteration (by means of Jacob's ladder $\vp_1(t)$, see \cite{3}) of the basic segment 
\bdis 
[\pi L,\pi L+U]=[\overset{0}{\wideparen{\pi L}},\overset{0}{\wideparen{\pi L+U}}], 
\edis  
and, finally, 
\be \label{1.4} 
\begin{split}
& \tilde{Z}^2(t)=\frac{{\rm d}\vp_1(t)}{{\rm d}t}=\frac{|\zf|^2}{\omega(t)}, \\ 
& \omega(t)=\left\{1+\mcal{O}\left(\frac{\ln\ln t}{\ln t}\right)\right\}\ln t, 
\end{split} 
\ee  
(see \cite{2}, (6.1), (6.7), (7.7), (7.8), (9.1)). 

\begin{remark} 
The components of the main $\zeta$-disconnected set 
\bdis 
\Delta(\pi L,U,1)=[\pi L,\pi L+U]\bigcup [\overset{1}{\wideparen{\pi L}},\overset{1}{\wideparen{\pi L+U}}]
\edis  
(for our case) are separated each from other by the gigantic distance $\rho$ (see \cite{3}, (5.12), comp. \cite{7}, (2.2)--(2.9)): 
\be\label{1.5} 
\rho\{[\pi L,\pi L+U]; [\overset{1}{\wideparen{\pi L}},\overset{1}{\wideparen{\pi L+U}}]\}\sim\pi(1-c)\frac{L}{\ln L},\ L\to\infty, 
\ee 
where $c$ stands for the Euler's constant. 
\end{remark}   

Further, let us remind that in our papers \cite{8} -- \cite{10} we have used an asymptotic form of corresponding exact complete hybrid formula (\ref{1.2}). In this paper we shall use directly the exact complete hybrid formula (\ref{1.2}) for our purposes. 

\subsection{} 

First we shall consider the following four sets 
\be \label{1.6} 
\begin{split}
& \{\zeta(ns)\},\ \{\Gamma(ns)\},\ \{\cn(ns,k)\},\ \{J_p(ns)\}, \\ 
& n\in\mbb{N},\ ns\in\mbb{C}\setminus \{ N,P\},\ k^2\in (0,1),\ p\in\mbb{Z}, 
\end{split} 
\ee  
where 
\bdis 
\{ N,P\}
\edis  
is the set of all zeros and poles of functions 
\bdis 
\zeta(s),\ \Gamma(s),\ \cn(s,k),\ J_p(s). 
\edis  
We obtain (for example) the following result in this direction: There are the sets 
\bdis 
{\overset{0}{\Omega}}^n_l,\ l=1,\dots,4,\ n\in\mbb{N} 
\edis  
such that we have the following infinite set  
\be \label{1.7} 
\begin{split}
& |\zeta(m{\overset{0}{{}s{}}}^m_1)| |\Gamma(m{\overset{0}{s}}^m_2)| |\cn(n{\overset{0}{s}}^n_3,k)|+|J_p(n{\overset{0}{s}}^n_4)| |\cn(m{\overset{0}{s}}^m_3,k)|= \\ 
& = |\zeta(n{\overset{0}{{}s{}}}^n_1)| |\Gamma(n{\overset{0}{s}}^n_2)| |\cn(m{\overset{0}{s}}^m_3,k)|+|J_p(m{\overset{0}{s}}^m_4)| |\cn(n{\overset{0}{s}}^m_n,k)|, \\ 
& (m,n)\in \mbb{N}^2,\ {\overset{0}{s}}^m_l\in{\overset{0}{\Omega}}^m_l,\ l=1,2,3,4 
\end{split}
\ee  
of exact meta-functional equations as new species generated by the mother formula (\ref{1.2}). Next, it follows from (\ref{1.7})  that on the infinite set 
\be \label{1.8} 
\begin{split}
& \{K(m,n)\}=\\ 
& = \{|\zeta(m{\overset{0}{{}s{}}}^m_1)||\Gamma(m{\overset{0}{{}s{}}}^m_2)||\cn(n{\overset{0}{{}s{}}}^n_3,k)|+|J_p(n{\overset{0}{{}s{}}}^n_4)||\cn(m{\overset{0}{{}s{}}}^m_3,k)|\}, \\ 
& (m,n)\in\mbb{N}^2 
\end{split}
\ee  
with quite complicated elements, still the commutative law 
\be \label{1.9} 
K(m,n)=K(n,m) 
\ee  
holds true. 

\subsection{} 

Further, we shall consider  even more complicated case of the set of four tuples: 
\be \label{1.10} 
\begin{split}
& \{\zeta(s),\Gamma(s),\cn(s,k),J_p(s)\}, \\ 
& \{\Gamma(2s),\cn(2s,k),J_p(2s),\zeta(2s)\}, \\ 
& \{\cn(3s,k),J_p(3s),\zeta(3s),\Gamma(3s)\}, \\ 
& \{J_p(4s),\zeta(4s),\Gamma(4s),\cn(4s,k)\}, \\ 
& --- \\ 
& \{\zeta(5s),\Gamma(5s),\cn(5s,k),J_p(5s)\}, \\ 
& \{\Gamma(6s),\cn(6s,k),J_p(6s),\zeta(6s)\}, \\ 
& \{\cn(7s,k),J_p(7s),\zeta(7s),\Gamma(7s)\}, \\ 
& \{J_p(8s),\zeta(8s),\Gamma(8s),\cn(8s,k)\}, \\ 
& --- \\ 
& \{\zeta(9s),\Gamma(9s),\cn(9s,k),J_p(9s)\}, \\ 
&  \vdots  
\end{split}
\ee  
generated by the cyclical changes in order of symbols 
\bdis 
\zeta,\Gamma,\cn(k),J_p. 
\edis   
In this case there are sets (for example) 
\bdis 
\Omega_l^1,\Omega_l^4,\ l=1,2,3,4 
\edis  
of elements 
\bdis 
s_l^1,s_l^4 
\edis  
such that the following equation 
\be \label{1.11} 
\begin{split}
& |\zeta(s_1^1)||\Gamma(s_2^1)||\Gamma(4s_3^4)|+|\cn(s_3^1,k)||\cn(4s_4^4,k)|= \\ 
& = |\zeta(4s_2^4)||J_p(4s_1^4)||\cn(s_3^1,k)|+|J_p(s_4^1)||\Gamma(4s_3^4)|
\end{split}
\ee 
holds true. 

\begin{remark}
Of course, the exact meta-functional equation (\ref{1.11}) represents only one element of an infinite sets of corresponding descendans (new species) generated by the mother formula (\ref{1.2}) as we shall see. 
\end{remark} 

\begin{remark}
The symmetry similar to (\ref{1.9}) is no longer valid for (\ref{1.11}). 
\end{remark} 

\begin{remark}
However, there are also meta-functional equations of the type: 
\be \label{1.12} 
\begin{split}
& |\zeta(1s_1^1)||\Gamma(1s_2^1)||\cn(5s_3^5,k)|+|J_p(5s_4^5)||\cn(1s_3^1,k)|= \\ 
& = |\zeta(5s_1^5)||\Gamma(5s_2^5)||\cn(1s_3^1,k)|+|J_p(1s_4^1)||\cn(5s_3^5)|
\end{split} 
\ee 
and for them (\ref{1.9})-type symmetries hold true w.r.t. transposition $1\leftrightarrow 5$. 
\end{remark} 

\subsection{} 

Let us notice that we present in this paper a kind of general method of generating infinite sets of meta-functional equations as the set of descendants of the mother formula (\ref{1.2}). 

\begin{remark}
Genericity of the presented method is based on: 
\begin{itemize}
	\item[(A)]  infinity of a set of possibilities how to make a choice of the mother formula 
	\item[(B)] applicability of this method to infinite set of admissible infinite sets of types (\ref{1.6}), (\ref{1.10}), \dots 
\end{itemize}
\end{remark} 

Next, we give the list of stages leading us to new infinite sets of meta-functional equations in present paper: 
\begin{itemize}
	\item[(a)] We assign the set of two factorization formulas (see \cite{4} -- \cite{7}) to the set (\ref{1.1}), 
	\item[(b)] crossbreeding on the last set gives us the exact complete hybrid formula (\ref{1.2}) defined od the critical line $\sigma=\frac 12$, 
	\item[(c)] next, we assign corresponding sets of level-curves to sets (\ref{1.6}), (\ref{1.10}), 
	\item[(d)] by means of mentioned level-curves we assign infinite sets of meta-functional equations to the mother formula (\ref{1.2}) in such a way that every equation contains an identical (=neutral) factor (with respect to operation (c)), 
	\item[(e)] finally, the crossbreeding on the last sets gives mentioned meta-functional equations on the complex plane as new species generated by the mother formula (\ref{1.2}). 
\end{itemize}

\begin{remark}
Also this paper is based on new notions and methods in the theory of the Riemann's zeta-function we have introduced in the series of 52 papers concerning Jacob's ladders. These can be found in arXiv [math.CA] starting with the paper \cite{1}. 
\end{remark} 

\section{Set of four tuples with a simple ordering of the symbols $\zeta,\Gamma,\cn(k),J_p$} 

\subsection{} 

Now, we assign corresponding level-curves to each set in (\ref{1.6}). Namely, we define four sets of level-curves 
\bdis 
\{{\overset{0}{\Omega}}_l^n\},\ l=1,2,3,4,\ {\overset{0}{\Omega}}_l^n\in\mbb{C}, 
\edis  
where 
\be \label{2.1} 
\begin{split}
& {\overset{0}{\Omega}}_1^n={\overset{0}{\Omega}}_1^n(\vec{S}_1,|\zeta(ns)|), \\ 
& {\overset{0}{\Omega}}_2^n={\overset{0}{\Omega}}_2^n(\vec{S}_1,|\Gamma(ns)|), \\ 
& {\overset{0}{\Omega}}_3^n={\overset{0}{\Omega}}_3^n(\vec{S}_2,|\cn(ns,k)|), \\ 
& {\overset{0}{\Omega}}_4^n={\overset{0}{\Omega}}_4^n(\vec{S}_3,|J_p(ns)|), 
\end{split}
\ee  
and 
\be \label{2.2} 
\begin{split}
& \vec{S}_1=(U,\pi L;f_1,|\zeta_{0,5}|^2), \\ 
& \vec{S}_2=(U,\pi L;f_2,|\zeta_{0,5}|^2), \\ 
& \vec{S}_3=(U,\pi L;|\zeta_{0,5}|^2), 
\end{split}
\ee 
as the loci 
\be \label{2.3} 
\begin{split}
& |\zeta(n{\overset{0}{{}s{}}}_1^n)|=\tilde{Z}^2(\alpha_1^{1,1})=c_1, \\ 
& |\Gamma(n{\overset{0}{{}s{}}}_2^n)|=\sin^2\alpha_0^{1,1}=c_2, \\ 
& |\cn(n{\overset{0}{{}s{}}}_3^n,k)|=\cos^2\alpha_0^{2,1}=c_3, \\ 
& |J_p(n{\overset{0}{{}s{}}}_4^n)|=\tilde{Z}^2(\beta_1^{1})=c_4, 
\end{split}
\ee  
and 
\be \label{2.4} 
0<c_1\dots, c_4<+\infty 
\ee  
for every admissible (see (\ref{1.3})) and fixed $U,L,k,p$, where 
\be \label{2.5} 
{\overset{0}{{}s{}}}_l^n\in {\overset{0}{\Omega}}_l^n
\ee 

\begin{remark}
Let us compare conditions (\ref{2.4}) with corresponding inclusions in (\ref{1.6}). 
\end{remark} 

Consequently, we have the following (see (\ref{1.2}), (\ref{2.3}), (\ref{2.5})): 

\begin{mydef51}
There is the set 
\be \label{2.6} 
|\zeta(n{\overset{0}{{}s{}}}_1^n)||\Gamma(n{\overset{0}{{}s{}}}_2^n)|+\tilde{Z}^2(\alpha_1^{2,1})|\zeta(n{\overset{0}{{}s{}}}_3^n,k)|=|J_p(n{\overset{0}{{}s{}}}_4^n)|,\ n\in\mbb{N} 
\ee  
of exact meta-functional equations with the neutral factor 
\bdis 
\tilde{Z}^2(\alpha_1^{2,1}) 
\edis  
(comp. the point (d) in the subsection 1.4 of this text). 
\end{mydef51}  

\subsection{} 

Next, we can write (see (\ref{2.6}))  
\be \label{2.7} 
\begin{split}
& |\zeta(m{\overset{0}{{}s{}}}_1^m)||\Gamma(m{\overset{0}{{}s{}}}_2^m)|+\tilde{Z}^2(\alpha_1^{2,1})|\cn(m{\overset{0}{{}s{}}}_3^m,k)|=|J_p(m{\overset{0}{{}s{}}}_4^m)|, \\ 
& |\zeta(n{\overset{0}{{}s{}}}_1^n)||\Gamma(n{\overset{0}{{}s{}}}_2^n)|+\tilde{Z}^2(\alpha_1^{2,1})|\cn(n{\overset{0}{{}s{}}}_3^n,k)|=|J_p(n{\overset{0}{{}s{}}}_4^n)|, \\ 
& (m,n)\in\mbb{N}^2,\ m\not=n. 
\end{split}
\ee 

Now we use the operation of crossbreeding to obtain the following Theorem (comp. \cite{6}, \cite{7} on the set (2.7), in this context elimination of $\tilde{Z}^2(\alpha_1^{1,2})$) 

\begin{mydef11} 
	There is infinite set 
\be \label{2.8} 
\begin{split}
& |\zeta(m{\overset{0}{{}s{}}}_1^m)||\Gamma(m{\overset{0}{{}s{}}}_2^m)||\cn(n{\overset{0}{{}s{}}}_3^n,k)|+|J_p(n{\overset{0}{{}s{}}}_4^n)||\cn(m{\overset{0}{{}s{}}}_3^m,k)|=\\ 
& = 
|\zeta(n{\overset{0}{{}s{}}}_1^n)||\Gamma(n{\overset{0}{{}s{}}}_2^n)||\cn(m{\overset{0}{{}s{}}}_3^m,k)|+|J_p(m{\overset{0}{{}s{}}}_4^m)||\cn(n{\overset{0}{{}s{}}}_3^n,k)|, \\ 
& (m,n)\in\mbb{N}^2,\ m\not=n,\ k^2\in (0,1), p\in\mbb{Z}, 
\end{split}
\ee  
(for every admissible and fixed $k$ and $p$) of exact meta-functional equations as new species generated by the mother formula (\ref{1.2}). 
\end{mydef11} 

\subsection{} 

We see immediately that for elements of the set 
\be \label{2.9} 
\begin{split}
& \{K(m,n)\}=\{|\zeta(m{\overset{0}{{}s{}}}_1^m)||\Gamma(m{\overset{0}{{}s{}}}_2^m)||\cn(n{\overset{0}{{}s{}}}_3^n,k)|+|J_p(n{\overset{0}{{}s{}}}_4^n)||\cn(m{\overset{0}{{}s{}}}_3^m,k)|\}, \\ 
& (m,n)\in\mbb{N}^2 
\end{split}
\ee  
the following is true: 

\begin{mydef41} 
\be \label{2.10} 
K(m,n)=K(n,m),\ \forall\- (m,n)\in\mbb{N}^2, 
\ee  
(the case $m=n$ is the trivial one). 
\end{mydef41} 

\section{Set of four tuples with cyclically ordered symbols $\zeta,\Gamma,\cn(k),J_p$}

\subsection{} 

Further, we assign corresponding level curves to elements of the first four sets in (\ref{1.10}). Namely 
\bdis 
\Omega_l^m,\ m,l=1,2,3,4 
\edis  
as the loci 
\begin{align} \label{3.1} 
& \Omega_l^1: & & \Omega_l^2: \nonumber \\ 
& |\zeta(s_1^1)|=\tilde{Z}^2(\alpha_1^{1,1})=c_1, & & |\Gamma(2s_1^2)|=c_1, \nonumber  \\ 
& |\Gamma(s_2^1)|=\sin^2(\alpha_0^{1,1})=c_2, & & |\cn(2s_2^2,k)|=c_2, \nonumber \\ 
& |\cn(s_3^1,k)|=\cos^2(\alpha_0^{2,1})=c_3, & & |J_p(2s_3^2)|=c_3, \nonumber \\ 
& |J_p(s_4^1)|=\tilde{Z}^2(\beta_1^1)=c_4, & & |\zeta(2s_4^2)|=c_4, \nonumber \\ 
& s_l^1\in\Omega_l^1 & & s_l^2\in\Omega_l^2, \\ 
& \Omega_l^3: & & \Omega_l^4: \nonumber \\ 
& |\cn(3s_1^3,k)|=c_1, & & |J_p(4s_1^4)|=c_1, \nonumber \\ 
& |J_p(3s_2^3,k)|=c_2, & & |\zeta(4s_2^4)|=c_2, \nonumber \\ 
& |\zeta(3s_3^3)|=c_3, & & |\Gamma(4s_3^4)|=c_3, \nonumber \\ 
& |\Gamma(3s_4^3)|=c_4, & & |\cn(4s_4^4,k)|=c_4, \nonumber \\ 
& s_l^3\in\Omega_l^3 & & s_l^4\in\Omega_l^4, \nonumber  \\ 
& 0<c_1,c_2,c_3,c_4<+\infty, \nonumber 
\end{align} 
where (comp. (\ref{2.2})) 
\be \label{3.2} 
\Omega_1^1=\Omega_1^1(\vec{S}_1,|\zeta(1s)|),\dots,\Omega_4^4=\Omega_4^4(\vec{S}_3,|\cn(4s,k)|)
\ee  
for every admissible (see (\ref{1.3}), (\ref{2.8})) and fixed $U,L,k,p$. 

Thus, we have (see (\ref{1.2}), (\ref{3.1})) the following 

\begin{mydef52}
\be \label{3.3} 
|\zeta(s_1^1)||\Gamma(s_2^1)|+\tilde{Z}^2(\alpha_1^{2,1})|\cn(s_3^1,k)|=|J_p(s_4^1)|, 
\ee 
\be \label{3.4} 
|\Gamma(2s_1^2)||\cn(2s_2^2,k)|+\tilde{Z}^2(\alpha_1^{2,1})|J_p(2s_3^2,k)|=|\zeta(2s_4^2)|, 
\ee 
\be \label{3.5} 
|\cn(3s_1^3,k)||J_p(3s_2^3)|+\tilde{Z}^2(\alpha_1^{2,1})|\zeta(3s_3^3)|=|\Gamma(3s_4^3)|, 
\ee  
\be \label{3.6} 
|J_p(4s_1^4)||\zeta(4s_2^4)|+\tilde{Z}^2(\alpha_1^{2,1})|\Gamma(4s_3^4)|=|\cn(4s_4^4,k)|  
\ee  
of transmutations of the mother formula (\ref{1.2}) with the neutral factor $\tilde{Z}^2(\alpha_1^{2,1})$.
\end{mydef52}  

\subsection{} 

Now, we use the operation of crossbreeding (see \cite{6}, \cite{7}) on every two different elements of the set 
\be \label{3.7}  
\{(3.3), (3.4), (3.5), (3.6)\}. 
\ee  

\begin{remark}
Let the symbol 
\bdis 
(3.3)\times (3.6) \ \Rightarrow 
\edis 
stand for \emph{we obtain by crossbreeding of the elements (3.3) and (3.6)}. 
\end{remark} 

Consequently, we obtain the following statement. 

\begin{mydef12} 
There is the following set of exact meta-functional equations as other transmutations of the mother formula (\ref{1.2}): 
\be \label{3.8} 
\begin{split}
& (3.3)\times (3.4) \ \Rightarrow \\ 
& |\zeta(s_1^1)||\Gamma(s_2^1)||J_p(2s_3^2)|+|\zeta(2s_4^2)||\cn(s_3^1,k)|=\\ 
& = |\Gamma(2s_1^2)||\cn(2s_2^2,k)||\cn(s_3^1,k)|+|J_p(s_4^1)||J_p(2s_3^2)|, 
\end{split}
\ee  
\be \label{3.9} 
\begin{split}
	& (3.3)\times (3.5) \ \Rightarrow \\ 
	& |\zeta(s_1^1)||\zeta(3s_3^3)||\Gamma(s_2^1)|+|\Gamma(3s_4^3)||\cn(s_3^1,k)|=\\ 
	& = |\cn(s_3^1,k)||\cn(3s_1^3,k)||J_p(3s_2^3)|+|J_p(s_4^1)||\zeta(3s_3^3)|, 
\end{split}
\ee  
\be \label{3.10} 
\begin{split}
	& (3.3)\times (3.6) \ \Rightarrow \\ 
	& |\zeta(s_1^1)||\Gamma(s_2^1)||\Gamma(4s_3^4)|+|\cn(s_3^1,k)||\cn(4s_4^4,k)|=\\ 
	& = |\zeta(4s_2^4)||J_p(4s_1^4)||\cn(s_3^1,k)|+|J_p(s_4^1)||\Gamma(4s_3^4)|, 
\end{split}
\ee  
\be \label{3.11} 
\begin{split}
	& (3.4)\times (3.5) \ \Rightarrow \\ 
	& |\zeta(3s_3^3)||\Gamma(2s_1^2)||\cn(2s_2^2,k)|+|\Gamma(3s_3^4)||J_p(2s_3^2)|=\\ 
	& = |\cn(3s_1^3,k)||J_p(2s_3^2)||J_p(3s_2^3)|+|\zeta(3s_3^3)||\zeta(2s_4^2)|, 
\end{split}
\ee 
\be \label{3.12} 
\begin{split}
	& (3.4)\times (3.6) \ \Rightarrow \\ 
	& |\Gamma(2s_1^2)||\Gamma(4s_3^4)||\cn(2s_2^2,k)|+|\cn(4s_4^4,k)||J_p(2s_3^2)|=\\ 
	& = |J_p(2s_3^2)||J_p(4s_1^4)||\zeta(4s_2^4)|+|\zeta(2s_4^2)||\zeta(4s_3^4)|, 
\end{split}
\ee 
\be \label{3.13} 
\begin{split}
	& (3.5)\times (3.6) \ \Rightarrow \\ 
	& |\cn(3s_1^3,k)||J_p(3s_2^3)||\Gamma(4s_3^4)|+|\cn(4s_4^4,k)||\zeta(3s_3^3)|=\\ 
	& = |J_p(4s_1^4)||\zeta(3s_3^3)||\zeta(4s_2^4)|+|\Gamma(3s_3^4)||\Gamma(4s_3^4)|.  
\end{split}
\ee 
\end{mydef12}  

\subsection{} 

Cyclycal ordering of the symbols $\zeta,\Gamma,\cn(k),J_p$ in the set (\ref{1.10}) gives the following generalization of our Lemma 2. 

\begin{mydef53} 
\be \label{3.14} 
\begin{split}
& |\zeta[(4m+1)s_1^{4m+1}]||\Gamma[(4m+1)s_2^{4m+1}]|+\tilde{Z}^2(\alpha_1^{2,1})|\cn[(4m+1)s_3^{4m+1},k]| \\ 
& =|J_p[(4m+1),s_4^{4m+1}]|, \\ 
& |\Gamma[(4m+2)s_1^{4m+2}]||\cn[(4m+2)s_2^{4m+2},k]|+\tilde{Z}^2(\alpha_1^{2,1})|J_p[(4m+2)s_3^{4m+2}]| \\ 
& =|\zeta[(4m+2),s_4^{4m+2}]|, \\ 
& |\cn[(4m+3)s_1^{4m+3},k]||J_p[(4m+3)s_2^{4m+3}]|+\tilde{Z}^2(\alpha_1^{2,1})|\zeta[(4m+3)s_3^{4m+3}]| \\ 
& =|\Gamma[(4m+3),s_4^{4m+3}]|, \\ 
& |J_p[(4m+4)s_1^{4m+4},k]||\zeta[(4m+4)s_2^{4m+4}]|+\tilde{Z}^2(\alpha_1^{2,1})|\Gamma[(4m+4)s_3^{4m+4}]| \\ 
& =|\cn[(4m+4),s_4^{4m+4},k]|, \\ 
& m\in\mbb{N}_0, 
\end{split}
\ee  
where the symbols 
\bdis 
s_1^{4m+1}\in\Omega_1^{4m+1},\dots,s_4^{4m+1}\in\Omega_1^{4m+1},\dots 
\edis  
one obtains as a natural continuation of the scheme (\ref{3.1}), see (\ref{1.10}). 
\end{mydef53} 

Next, we use the symbol 
\bdis 
(4m+q)\times (4m+r) \ \Rightarrow \ ,\ m,n\in\mbb{N}_0, \ q,r=1,2,3,4 
\edis  
in the sense of our Remark 8, where $4m+q$ and $4m+r$ are sequential numbers of equations in (\ref{3.14}). Consequently, we obtain the following generalization of our Theorem 2. 

\begin{mydef13} 
There are six infinite classes of exact meta-functional equations as new species generated by the mother formula (\ref{1.2}): 
\be \label{3.15} 
\begin{split}
& (4m+1)\times(4n+2) \ \Rightarrow \\ 
& |\zeta[(4m+1)s_1^{4m+1}]||\Gamma[(4m+1)s_2^{4m+1}]||J_p[(4n+2)s_3^{4n+2}]|+ \\ 
& + |\zeta[(4n+2)s_4^{4n+2}]||\cn[(4m+1)s_3^{4m+1},k]|= \\ 
& = |\Gamma[(4n+2)s_1^{4n+2}]||\cn[(4n+2)s_2^{4n+2},k]||\cn[(4m+1)s_3^{4m+1},k]|+ \\ 
& + |J_p[(4m+1)s_4^{4m+1}]||J_p[(4n+2)s_3^{4n+2}]|, \\ 
& \vdots \\ 
& (4m+3)\times(4n+4) \ \Rightarrow \\ 
& |\cn[(4m+3)s_1^{4m+3},k]||J_p[(4m+3)s_2^{4m+3}]||\Gamma[(4n+4)s_2^{4n+4}]|+ \\ 
& + |\cn[(4n+4)s_4^{4n+4},k]||\zeta[(4m+3)s_3^{4m+3}]|= \\ 
& = |J_p[(4n+4)s_1^{4n+4}]||\zeta[(4m+3)s_3^{4m+3}]||\zeta[(4n+4)s_2^{4n+4}]|+ \\ 
& + |\Gamma[(4m+3)s_4^{4m+3}]||\Gamma[(4n+4)s_3^{4n+4}]|, \\ 
& m,n\in\mbb{N}_0. 
\end{split}
\ee 	
\end{mydef13}  

\subsection{} 

Further, as a subset of such four consecutive sets in (\ref{1.10}) that corresponds to the rows 
\bdis 
4m+1,4m+2,4m+3,4m+4 
\edis  
we call the $m$-th cell of the set (\ref{1.10}). 

\begin{remark}
We may assume that results of some interactions between the $m$-th and $n$-th cells are expressed by the equations (\ref{3.15}). Namely: 
\begin{itemize}
	\item[(a)] If $m=n$, then we have the case of internal interactions between all mutually different elements of the $m$-th cell. 
	\item[(b)] If $m\not=n$, then we have the case of external interactions between every element (=corresponding four tuple) of the $m$-th cell with every element of the $n$-th cell. 
\end{itemize}
\end{remark} 

\subsection{} 

Next, we present some property connected with the external interactions mentioned above. Namely, we have (similarly to (\ref{3.15})) 
\be \label{3.16} 
\begin{split}
& (4m+1)\times (4n+1) \ \Rightarrow \\ 
& |\zeta[(4m+1)s_1^{4m+1}]||\Gamma[(4m+1)s_2^{4m+1}]||\cn[(4n+1)s_3^{4n+1},k]|+ \\ 
& + |J_p[(4n+1)s_4^{4n+1}]||\cn[(4m+1)s_3^{4m+1},k]|= \\ 
& = |\zeta[(4n+1)s_1^{4n+1}]||\Gamma[(4n+1)s_2^{4n+1}]||\cn[(4m+1)s_3^{4m+1},k]|+ \\ 
& + |J_p[(4m+1)s_4^{4m+1}]||\cn[(4n+1)s_2^{4n+1},k]|.
\end{split}
\ee  
Now, if we put 
\bdis 
\begin{split}
& G(4m+1,4n+1)= \\ 
& |\zeta[(4m+1)s_1^{4m+1}]||\Gamma[(4m+1)s_2^{4m+1}]||\cn[(4n+1)s_3^{4n+1},k]|+ \\ 
& + |J_p[(4n+1)s_4^{4n+1}]||\cn[(4m+1)s_3^{4m+1},k]|, 
\end{split}
\edis 
then we obtain for the set 
\bdis 
\{G(4m+1,4n+1)\},\ m,n\in\mbb{N}_0,\ m\not=n 
\edis  
the following 
\begin{mydef42} 
\be \label{3.17} 
\begin{split}
& G(4m+1,4n+1)=G(4n+1,4m+1), 
\end{split}
\ee  
and similarly 
\be \label{3.18} 
\begin{split}
	& G(4m+q,4n+q)=G(4n+q,4m+q),\ q=2,3,4.  
\end{split}
\ee  
\end{mydef42}

I would like to thank Michal Demetrian for his moral support of my study of Jacob's ladders.

\end{document}